\documentclass{ifacconf}

\usepackage{amssymb}
\newcommand*{\qedb}{\hfill\ensuremath{\square}}%

\usepackage{amsmath}
\usepackage{color}
\usepackage[shortlabels]{enumitem}

\usepackage{balance}
\usepackage{tikz}
\usepackage{tkz-graph}

\usepackage{graphicx}      % include this line if your document contains figures
\usepackage{natbib}        % required for bibliography
\begin{document}
\begin{frontmatter}

\title{Distributed Linear Quadratic Tracking Control for Leader-Follower Multi-Agent Systems: A Suboptimality Approach %\thanksref{footnoteinfo}
} 

%\thanks[footnoteinfo]{The work of Junjie Jiao was supported by China Scholarship Council (CSC).}

\author[First]{Junjie Jiao} 
\author[First]{Harry L. Trentelman} 
\author[First]{M. Kanat Camlibel}

\address[First]{Bernoulli Institute for Mathematics, Computer Science and Artificial Intelligence, University of Groningen, Groningen, 9700 AV, The Netherlands (e-mail: \{ j.jiao, h.l.trentelman, m.k.camlibel\}@rug.nl).}

\begin{abstract}                % Abstract of not more than 250 words.
In this paper, we extend the results from \cite{8736845} on  distributed linear quadratic control for leaderless multi-agent systems to the case of distributed linear quadratic  tracking control for leader-follower multi-agent systems.  
Given one autonomous leader and a number of homogeneous followers, we introduce an associated global quadratic cost functional. 
We assume that the leader shares its state information with at least one of the followers and the communication between the followers is represented by a connected simple undirected graph.
Our objective is to design distributed control laws such that the controlled network reaches tracking consensus  and, moreover, the associated cost is smaller than a given tolerance for all initial states bounded in norm by a given radius.
We establish a centralized design method for computing such suboptimal control laws, involving the solution of a single Riccati inequality of dimension equal to the dimension of the local agent dynamics, and the smallest and the largest eigenvalue of a given positive definite matrix involving the underlying graph.
The proposed design method is illustrated by a simulation example.
\end{abstract}


\end{frontmatter}
%===============================================================================
\section{Introduction}
Distributed control for multi-agent systems has drawn much attention in the past two decades due to its practical applications, e.g., formation control, intelligent transportation systems and power grids.
In the literature, basically two types of multi-agent systems are considered, namely leaderless multi-agent systems and leader-follower multi-agent systems. 
In the leaderless case, the local agents reach agreement  which depends on the dynamics of all agents  (\cite{Olfati-Saber2004}, \cite{harry_2013}). In the leader-follower case,  the states or the outputs of  the followers  track that of the leader   (\cite{HONG20061177}, \cite{NI2010209}).
One of the attractive directions in distributed control for multi-agent systems is to design distributed control laws that minimize certain global or local performances, while reaching an agreement for the controlled network.

In the past, quite some work has been devoted to distributed linear quadratic (LQ) optimal control for leaderless multi-agent systems.
In  \cite{tuna2008lqr}, an LQR based method was used to design distributed synchronizing control laws for a multi-agent system, without taking any performance into consideration.
In \cite{tamas2008}, suboptimal distributed stabilizing control laws were established for a multi-agent system with general agent dynamics with respect to an associated global cost functional, while in \cite{wei_ren2010}, the optimal synchronizing control gain was computed for leaderless multi-agent systems with single integrator agent dynamics.
In the meantime, the distributed LQ control problem was also considered in \cite{Semsar2009} by utilizing a game theoretic approach, in \cite{kristian2014} by adopting an inverse optimal approach, and later  in  \cite{8736845} by employing a suboptimality approach. 
For other papers related to this topic,  see also  \cite{huaguang_zhang2015} and \cite{Jiao2019}.

On the other hand, distributed LQ tracking control for leader-follower multi-agent systems has also attracted much attention.
In \cite{hongwei_zhang2011}, distributed synchronizing control laws were established using an LQR based approach without optimizing any performance.
Later on, in \cite{Guaranteed_LF2015}, suboptimal distributed control laws were proposed for achieving   guaranteed cost.
In \cite{Nguyen2015}, a hierarchical LQR based method was presented  to design suboptimal synchronizing control laws for leader-follower systems, and an inverse optimal approach was introduced in \cite{kristian2014}, see also \cite{huaguang_zhang2015}.

In the present paper we extend the results from \cite{8736845} on   distributed LQ control  for leaderless multi-agent systems to the case of distributed LQ tracking control for leader-follower multi-agent systems. Given a leader-follower system with one autonomous leader and a number of followers, we introduce an associated global quadratic cost functional. 
We assume that the leader shares its state information with at least one of the followers, and the communication between the followers is represented by a  connected simple undirected graph.
Our aim is then to  design  distributed diffusive control laws such that the controlled network reaches tracking consensus, i.e., the states of the followers track the state of the leader asymptotically and the associated cost is smaller than an a priori given upper bound.

%%%%%%%%%%%%%%%%%%T%%%%%%%%%%%%

The outline of this paper is as follows. Section \ref{sec_preliminary} provides
some  preliminaries on graph theory and quadratic performance analysis for  linear systems.
In Section \ref{sec_problem}, we formulate
the suboptimal distributed linear quadratic tracking control problem for leader-follower multi-agent systems. 
We then address this  suboptimal distributed tracking control problem in Section \ref{sec_mas}. A simulation example is presented in Section \ref{sec_simulation} to illustrate our design method.
Finally,  Section \ref{sec_conclusion} concludes this paper.

\subsection*{Notation}
We denote by $\mathbb{R}$ the field of real numbers, and by $\mathbb{R}^n$ the  $n$-dimensional  real Euclidean space.
The column vector $\mathbf{1}_n\in \mathbb{R}^n$ denotes the vector whose entries are all equal to $1$.
For $x \in \mathbb{R}^n$, we define  its Euclidean norm $\|x\| := \sqrt{x^\top x}$. 
For a given $r >0$, we denote by $B(r) := \{x \in \mathbb{R}^n \mid \|x\| \leq r\}$ the closed ball of radius $r$.
We denote by $\mathbb{R}^{n\times m}$ the space of  real $n\times m$ matrices. 
For a given matrix $A$, its transpose and inverse (if it exists) are denoted by $A^{\top}$ and $A^{-1}$, respectively.
We denote by $I_n$ the identity matrix of dimension $n \times n$. 
The Kronecker product of two matrices $A$ and $B$ is denoted by $A\otimes B$, which has the property that $(A_1\otimes B_1)(A_2\otimes B_2) = A_1 A_2\otimes B_1 B_2 $. 
For a given symmetric matrix $P$ we denote $P>0$ if it is positive definite and $P < 0$ if it is negative definite.
By  $\text{diag} (  a_1, a_2, \ldots, a_n )$,
%\todoiny{?}{I would write $\text{diag} (  a_1, a_2, \ldots, a_n )$. Curly brackets are used for sets. }
we denote the $n\times n$ diagonal matrix with $a_1, a_2, \ldots, a_n $ on the diagonal.

\section{Preliminaries}\label{sec_preliminary}
%In this section, we will review basic graph theory and the linear quadratic performance analysis of autonomous  systems.

\subsection{Graph Theory}
%\blue{
In this paper,  a (directed) graph is a tuple $\mathcal{G} = (\mathcal{V},\mathcal{E})$ with nonempty node set $\mathcal{V} = \{1, 2, \ldots, N \}$ and edge set $\mathcal{E} \subset \mathcal{V} \times \mathcal{V}$.
The edge from node $i$ to node $j$ is represented by the pair $(i, j) \in \mathcal{E}$ with $i, j\in \mathcal{V}$. 
We say  the graph is simple if the edge set $\mathcal{E}$ only contains edges of the form $(i,j)$ with $i \neq j$. 
The graph is called undirected if $(i, j) \in \mathcal{E}$ implies $(j, i) \in \mathcal{E}$. 
%
%Given node $i$, we denote its neighboring set by $\mathcal{N}_i :=\{ j \in \mathcal{V}  \mid (i, j)\in \mathcal{E} \}$.
%
The adjacency matrix of the graph $\mathcal{G}$ is defined as $\mathcal{A} = [a_{ij}]$ with $a_{ij} =1$ if there is an edge between the nodes $i$ and $j$, and $a_{ij} = 0$ otherwise.
For simple graphs, $a_{ii} = 0$ for all $i$. 
Furthermore, a graph $\mathcal{G}$ is undirected if and only if $\mathcal{A}$ is symmetric.
The Laplacian matrix is defined as ${L} = D - \mathcal{A}$, where  $D = \text{diag} ( d_1,d_2,\ldots, d_N )$ with $d_{i} = \sum_{j=1}^{N} a_{ij}$  the degree matrix of $\mathcal{G}$.
The Laplacian matrix $L$ of an undirected graph is symmetric and consequently only has real eigenvalues. Furthermore, all eigenvalues are nonnegative and $0$ is an eigenvalue of $L$.
The graph is connected if and only if $0$ is a simple eigenvalue of $L$. 
%
%In the sequel we will assume that $\mathcal{G}$ is connected simple undirected graph. 
%In that case the eigenvalues of $L$ can be ordered in increasing order as $0=\lambda_1 < \lambda_2 \leq \cdots \leq \lambda_N$ and there exists an orthogonal matrix $U$ such that 
%$U^{\top}LU = \text{diag} ( 0, \lambda_2, \ldots, \lambda_N )$.
%
%\todoiny{?}{I would write $\text{diag} (  a_1, a_2, \ldots, a_n )$. Curly brackets are used for sets. }
%Moreover, we have $U = \left( \frac{1}{\sqrt{N}}\mathbf{1}_N\quad U_2 \right)$ and $U_2 U^{\top}_2 = I_N - \frac{1}{N}\mathbf{1}_N\mathbf{1}_N^{\top}$.

For a connected simple undirected graph  $\mathcal{G}$, we review the following result:
\begin{lem}[\cite{HONG20061177}]\label{Gamma}
	Let $\mathcal{G}$ be a connected simple undirected graph with Laplacian matrix $L$. 
	Let $g_1,g_2,\ldots,g_N$ be non-negative real numbers with at least one $g_i >0$. Define $G = \textnormal{diag} (g_1,g_2,\ldots,g_N)$. 
	Then the matrix $L+G$ is positive definite.
\end{lem}

\subsection{Quadratic Performance of Linear Autonomous  Systems}
In this subsection, we analyze the quadratic performance  of a linear autonomous  system.
Consider the autonomous system 
\begin{equation}\label{sys_auto}
\dot{x}(t) =\bar{A}x(t),\quad x(0) = x_0
\end{equation}
where  $\bar{A}\in \mathbb{R}^{n\times n}$ and $x\in \mathbb{R}^{n}$ is the state.
We consider the quadratic performance of system (\ref{sys_auto}),  given by
\begin{equation}\label{cost_auto}
J = \int_{0}^{\infty} x^{\top}(t) \bar{Q} x(t) \ dt
\end{equation}
where  $\bar{Q}\geq 0$ is a given real weighting matrix.
Note that the performance $J$ is finite if  system (\ref{sys_auto}) is stable, i.e., $\bar{A}$ is Hurwitz.

%
%We are interested in finding conditions such that the performance (\ref{cost_auto}) of system (\ref{sys_auto}) is smaller than a given upper bound. 
%%
%For this, we have the following theorem:

The following well-known result (\cite{algebraic_1997} and \cite{8736845}) provides a {\em necessary} and {\em sufficient} condition such that, for a given tolerance $\gamma >0$, the  performance (\ref{cost_auto}) satisfies $J < \gamma$.

\begin{thm}\label{thm_autonomous}
	Consider  system (\ref{sys_auto}) with associated  performance (\ref{cost_auto}).
	For given $\gamma > 0$, we have that $\bar{A}$ is Hurwitz and $J < \gamma$ if and only if there exists  $P>0$ satisfying
	%	\todoiny{?}{Replace "to the inequalities" by "satisfying"}
	\begin{align}
	\bar{A}^{\top} P  + 	P\bar{A} + \bar{Q} &< 0,\label{lya_ineq_3} \\
	x_0^{\top} Px _0  &< \gamma. \label{conditions_p}
	\end{align}
\end{thm}

In the next section, we  formulate the problem that we will address in this paper.

\section{Problem Formulation}\label{sec_problem}
In this paper, we consider a leader-follower multi-agent system, consisting of one leader and  $N$  followers. 
The dynamics of the leader is represented by
the linear time-invariant autonomous system 
\begin{equation}\label{sys_leader}
\dot{x}_r(t) = Ax_r(t),\quad x_r(0) = x_{r0}.
\end{equation}
where $A\in \mathbb{R}^{n\times n}$, $x_r\in \mathbb{R}^{n}$ is the state of the leader and $x_{r0}$ is its initial state.
The dynamics of the followers are identical and  represented by the linear time-invariant systems
\begin{equation}\label{sys_follower}
\dot{x}_i(t) = Ax_i(t) + Bu_i(t),\quad x_i(0) = x_{i0}, \quad i=1, 2, \ldots,N
\end{equation}
where $A\in \mathbb{R}^{n\times n}$, $B\in \mathbb{R}^{n\times m}$, and $x_i\in \mathbb{R}^{n}, u_i \in \mathbb{R}^{m}$ are the state and input of follower $i$, respectively, and $x_{i0}$ is its initial state.
Throughout this paper, we assume  that the pair $(A,B)$ is stabilizable.
Moreover, we make the following two standard assumptions regarding the communication between the leader and the followers:
\begin{assum}\label{assum_leader}
	We assume that  at least one follower receives the state information of the leader.
\end{assum}

\begin{assum}\label{assum_follower}
	We also assume that the underlying graph $\mathcal{G}$ of the communication between the followers is a connected simple undirected graph.
\end{assum}

%
%Note that, although the followers have identical dynamics, we allow the initial states of the followers to differ. 
%

We consider the infinite horizon distributed linear quadratic tracking control problem for the leader-follower system \eqref{sys_leader} and \eqref{sys_follower}, where the global cost functional integrates the weighted quadratic difference of states between every follower and its neighbors and the weighted quadratic difference of states between the leader and the followers communicating with the leader, and where the cost functional also penalizes the inputs in a quadratic form.

Note that, as mentioned in Assumption \ref{assum_leader}, at least one follower receives the state information of the leader. 
Thus, the leader-follower system \eqref{sys_leader} and \eqref{sys_follower} can be interconnected by a distributed diffusive control law of the form 
\begin{equation}\label{control}
u_i (t)= K\sum_{j=1}^{N} a_{ij}(x_i(t) - x_j (t)) + K g_i (x_i(t) - x_r(t))
\end{equation} 
where  $a_{ij}$ is the $ij$-th entry of the adjacency matrix $\mathcal{A}$ of the underlying graph $\mathcal{G}$, $K \in \mathbb{R}^{m\times n}$ is an identical feedback gain for all followers and  we have $g_i > 0$ for at least one $i=1,2, \ldots, N$. 
%Note that, as mentioned in Assumption \ref{assum_leader}, at least one follower receives the state information of the leader, thus, 
%
Accordingly, the cost functional  considered in this paper is given by
\begin{equation}\label{cost}
\begin{aligned}
J(u) &= \int_{0}^{\infty}\frac{1}{2}\sum_{i=1}^{N}\sum_{j=1}^N a_{ij}(x_i-x_j)^{\top} Q (x_i-x_j) \\
&\qquad + \sum_{i=1}^{N}  g_{i}(x_i-x_r)^{\top} Q (x_i-x_r) + \sum_{i=1}^{N}u_i^{\top}R u_i \ dt
\end{aligned}
\end{equation}
where $Q\geq 0$ and $R > 0$ are given real weighting matrices of suitable dimensions.

The  distributed linear quadratic tracking problem is then the problem of minimizing  the cost functional (\ref{cost})  for all initial states $x_{r0}$ and  $x_{i0}$, $i=1,2,\ldots,N$ over all distributed diffusive control laws \eqref{control} such that the states of all followers track  the state of the leader asymptotically.
In that case we say the network reaches {\em tracking consensus}: 
\begin{defn}\label{def_consensus}
	We say the  control law \eqref{control} achieves tracking consensus for the leader-follower system \eqref{sys_leader} and \eqref{sys_follower} if for all $i =1,2,\ldots,N$ and for all initial states   $x_{r0}$ and $x_{i0}$, we have $$x_i(t)-x_r(t) \rightarrow 0
	\text{  as  } t \rightarrow \infty.$$
\end{defn}

Due to the distributed nature of the control law \eqref{control} as imposed by the network topology, the distributed linear quadratic tracking problem is a non-convex optimization problem (\cite{lunze_conf_2014}). 
It is therefore difficult, if not impossible, to find a closed form solution for an optimal controller, or such optimal controller may not even exist.
Therefore, in this paper we will  design distributed control laws which solve a {\em suboptimal} version of this problem. 

To proceed, for the $i$th follower we introduce the following error state
\begin{equation*}
e_i = x_i -x_r, 
\end{equation*}
for $i = 1, 2, \ldots,N.$
Subsequently, the dynamics of $e_i$ is given by
\begin{equation}\label{error_i}
\dot{e}_i	= A e_i + Bu_i,  \quad i = 1, 2, \ldots,N.
\end{equation}
Denoting $x = \left( x_1^{\top},\ldots,x_N^{\top} \right)^{\top}$, $u = \left( u_1^{\top},\ldots,u_N^{\top} \right)^{\top}$, and $e = \left( e_1^{\top},\ldots,e_N^{\top} \right)^{\top}$,
we can then rewrite the  error system \eqref{error_i} in compact form as
\begin{equation}\label{error}
\dot{e} = (I_N\otimes A) e + (I_N\otimes B)u,\quad e(0) =e_0.
\end{equation}
Note that 
\[
	e = x - \mathbf{1}_N\otimes x_r.
\]
Correspondingly, by using the fact $(L\otimes K)(\mathbf{1}\otimes x_r) =0,$ the control law \eqref{control} can be given by
\begin{equation}\label{control_error}
u(t) = (\Gamma \otimes K) e
\end{equation}
where $\Gamma = L+G$ and $G = \textnormal{diag}(g_1,g_2,\ldots, g_N)$.
%By substituting the control law \eqref{control_error} into the error system \eqref{error}, we have 
%\begin{equation}\label{error_closed}
%\dot{e} = (I_N\otimes A  + (L+G) \otimes BK ) e,\quad e(0) =e_0.
%\end{equation}
%
Similarly, the cost functional \eqref{cost} can be written in terms of $e$ and $u$ as
\begin{equation}\label{cost_error}
J(u) = \int_{0}^{\infty}  e^\top (\Gamma  \otimes Q) e + u^\top (I_N \otimes R) u ~ dt.
\end{equation}
Now, by substituting the control law \eqref{control_error} into the error dynamics \eqref{error}, we obtain the closed-loop error system
\begin{equation}\label{sys_error}
\dot{e} = (I_N\otimes A +  \Gamma \otimes BK)e,\quad e(0) =e_0.
\end{equation}
and the associated cost is now given by 
\begin{equation}\label{cost_error_K}
J(K) = \int_{0}^{\infty}  e^\top \Big(\Gamma  \otimes Q  +  \Gamma^2 \otimes K^\top R K \Big) e ~ dt
\end{equation}
Note that  the  controlled leader-follower system  \eqref{sys_leader} and \eqref{sys_follower} reaches tracking consensus, i.e., the states of all followers track the state of the leader asymptotically, if and only if the error dynamics \eqref{sys_error} is stable.

Let 
\begin{equation}\label{error_initial}
B(r) = \{e_0 \in \mathbb{R}^{nN} \mid \| e_0 \| \leq r\}
\end{equation}
be the closed ball of radius $r$ in the state space $\mathbb{R}^{nN}$ of the error system \eqref{sys_error}.
Then, for the leader-follower system \eqref{sys_leader} and \eqref{sys_follower} with initial states  such that the error initial state is contained in a closed ball of a given radius, we want to design a distributed diffusive controller such that tracking consensus is achieved and, for all initial states satisfying \eqref{error_initial}, the associated cost is smaller than an a priori given upper bound.
Thus, the problem that we will address is the following:
\begin{prob} \label{prob1}
	Consider the leader-follower multi-agent system \eqref{sys_leader} and \eqref{sys_follower} and the associated cost functional \eqref{cost}. Let $r >0$ be a given radius and let $\gamma >0$ be an a priori given upper bound for the cost.
	The problem is to find a distributed diffusive control law of the form \eqref{control} such that the controlled leader-follower system reaches tracking consensus and, for all initial conditions $x_0$ and $x_{r0}$ such that $e_0 = x_0 - x_{r0}$ satisfies \eqref{error_initial}, the associated cost \eqref{cost} is smaller than the given upper bound, i.e., $J(K) < \gamma$.
\end{prob}

%\begin{rem}
%Note that we could also have formulated the alternative problem of finding a suboptimal controller for a {\em single, given, initial state $x_0$}. In fact, this would be closer to the classical linear quadratic problem, which is usually formulated as the problem of minimizing the cost functional for a {\em given} initial state $x_0$. In that context, however, the optimal controller is a state feedback that turns out to be optimal {\em for all initial states}. In order to capture in our problem formulation this property of being optimal for all initial states, we have formulated Problem \ref{prob1} in terms of initial states contained in a ball of a given radius.
%\end{rem}

%
%

%
\section{Suboptimal Control Design for Leader-Follower Multi-Agent Systems}\label{sec_mas}
%

%As formulated in Problem \ref{prob1}, given a desired upper bound $\gamma >0$, for multi-agent system  \eqref{sys_leader} and \eqref{sys_follower} with initial relate states $e = \left( e_1^{\top},\ldots,e_N^{\top} \right)^{\top} $ between the followers and the leader contained in the closed ball $ e \in B(r)$ of given radius $r$ we want to design a control law of the form  \eqref{control_error} such that the error system \eqref{error} is internally stable  and, moreover, for all $ \in B(r)$ the associated cost  \eqref{cost_error_K} is smaller than the given upper bound, i.e., $J(K)<\gamma$.
In this section, we will address Problem \ref{prob1} and provide a  suitable control design method.
As  mentioned before, the distributed control law \eqref{control} achieves tracking consensus and suboptimal performance for the leader-follower system  \eqref{sys_leader} and \eqref{sys_follower}  with respect to the given tolerance on the cost functional \eqref{cost} if and only if the error dynamics \eqref{sys_error} is stable and $J(K) <\gamma$.

Now, let $U\in \mathbb{R}^{N\times N}$ be an orthogonal matrix that diagonalizes $\Gamma = L+G$. 
Define 
\[
U^{\top}\Gamma U :=  \Lambda = \text{diag} ( \lambda_1, \lambda_2,\ldots,\lambda_N ).
\]
It follows from Lemma \ref{Gamma} that  $\lambda_i>0$  for all $i=1,2,\ldots,N$.
%
%\todoiny{?}{I would write $\text{diag} (  a_1, a_2, \ldots, a_n )$. Curly brackets are used for sets. }
To simplify the problem formulated in the previous section, by applying the state transformation $\bar{e} =(U^{\top}\otimes I_n)e$,  system \eqref{sys_error} becomes
\begin{equation}\label{error_bar}
\dot{\bar{e}} = \left(I_N\otimes A + \Lambda \otimes BK\right) \bar{e},\quad \bar{e}(0) =\bar{e}_0 
\end{equation}
where $\bar{e} = \left( \bar{e}_1^{\top},\ldots,\bar{e}_N^{\top} \right)^{\top}$.
In terms of the transformed variable, the cost \eqref{cost_error_K} is then given by
\begin{equation}\label{cost_error_bar}
{J}(K) = \int_{0}^{\infty}\sum_{i=1}^{N} \bar{e}_i^{\top} (\lambda_i Q + \lambda_i^2 K^{\top} R K) \bar{e}_i\ dt.
\end{equation}
Note that the transformed states $\bar{e}_i$, $i = 1,2,\ldots, N$ appearing in system  \eqref{error_bar} and cost \eqref{cost_error_bar} are decoupled from each other.
Then we can write  system \eqref{error_bar} as
\begin{equation}\label{closed_loop}
\dot{\bar{e}}_i = (A +\lambda_i BK)\bar{e}_i ,\quad i = 1,2,\ldots,N. 
\end{equation}
Also, the cost \eqref{cost_error_bar} equals 
\begin{equation}\label{cost_cost}
{J}(K)= \sum_{i=1}^{N}{J}_i(K)
\end{equation}
with
\begin{equation}\label{coco}
{J}_i(K) = \int_{0}^{\infty} \bar{e}_i^{\top} (\lambda_i Q + \lambda_i^2 K^{\top}RK) \bar{e}_i\ dt, \quad i = 1, 2, \ldots,N.
\end{equation}
Clearly, the controlled leader-follower system \eqref{sys_leader} and \eqref{sys_follower}  reaches tracking consensus with control law \eqref{control} if and only if, for $i=1,2,\ldots,N$, the systems \eqref{closed_loop} are stable.
In addition, the control law \eqref{control} is suboptimal if $J(K) < \gamma$.

So far, we have transformed the problem of distributed suboptimal control for the leader-follower system \eqref{sys_leader} and \eqref{sys_follower} into the problem of finding one single static feedback gain  $K\in \mathbb{R}^{m\times n}$ such that the systems \eqref{closed_loop} are stable  for $i=1,2,\ldots,N$ and $J(K) < \gamma$. 
Since the pair $(A, B)$ is stabilizable, there exists such a feedback gain $K$ (\cite{zhongkui_li_unified_2010}, \cite{hongwei_zhang2011}).

The following lemma then provides a necessary and sufficient condition for a given  feedback gain $K$ to stabilize all systems \eqref{closed_loop}  and for given initial states guarantee that $J(K) < \gamma$.

\begin{lem}\label{lem_1}
	Let $K$ be a feedback gain.
	Consider the systems   \eqref{closed_loop} with given initial states $\bar{e}_{10}, \bar{e}_{20}, \ldots , \bar{e}_{N0}$ and associated cost functionals \eqref{cost_cost} and \eqref{coco}. 
	Let $\gamma>0$.
	Then all systems (\ref{closed_loop}) are stable and ${J}(K)<\gamma$ if and only if 
	there exist $P_i >0$ satisfying
	\begin{equation}
	(A +  \lambda_i B K)^{\top}P_i + P_i (A + \lambda_i B K)  + \lambda_i Q +\lambda_i^2 K^{\top}R K  <0 \label{are_ineq}
	\end{equation}
	and
	\begin{equation}
	\sum_{i=1}^{N} \bar{e}_{i0}^{\top} P_i \bar{e}_{i0}  <\gamma, \label{initial_condition_n_1}
	\end{equation}
	for $i= 1, 2, \ldots, N$, respectively. 
\end{lem}

\begin{pf}
	($\Leftarrow$)
	Since (\ref{initial_condition_n_1}) holds, there exist $\gamma_i := \bar{e}_{i0}^{\top} P_i \bar{e}_{i0} +\epsilon_i$   with  sufficiently small $\epsilon_i > 0$, $i=1,2,\ldots,N$ such that $\sum_{i=1}^{N} \gamma_i < \gamma$. 
	Because there exists $P_i>0$ such that (\ref{are_ineq}) and $\bar{e}_{i0}^{\top} P_i \bar{e}_{i0} < \gamma_i$ holds for all $i =1,2,\ldots,N$, by taking $\bar{A} = A + \lambda_i BK$ and $\bar{Q} = \lambda_i Q +\lambda_i^2 K^{\top}R K$, it follows from Theorem \ref{thm_autonomous} that all systems \eqref{closed_loop} are stable and $J_i(K) <\gamma_i$ for $i = 1,2,\ldots,N$. 
	Since $J(K) = \sum_{i=1}^{N} J_i(K)$, this implies that $J(K) < \sum_{i=1}^{N} \gamma_i < \gamma$.
	
	($\Rightarrow$)
	Since $J(K)<\gamma$ and $J(K) = \sum_{i=1}^{N} J_i(K)$, there exist $\gamma_i := J_i(K)+\epsilon_i$ with sufficiently small $\epsilon_i > 0$, $i=1,2,\ldots,N$ such that $\sum_{i=1}^{N} \gamma_i < \gamma$.
	Because all systems (\ref{closed_loop}) are stable and $J_i(K) <\gamma_i$ for $i=1,2,\ldots, N$, by taking $\bar{A} = A + \lambda_i BK$ and $\bar{Q} = \lambda_i Q +\lambda_i^2 K^{\top}R K$, it again follows from Theorem \ref{thm_autonomous} that there exist  $P_i>0$ such that (\ref{are_ineq})  and  $\bar{x}_{i0}^{\top} P_i \bar{x}_{i0} < \gamma_i$ hold for all $i =1,2,\ldots,N$. 
	Since $\sum_{i=1}^{N} \gamma_i < \gamma$, this implies that $\sum_{i=1}^{N} \bar{x}_{i0}^{\top} P_i \bar{x}_{i0}  <\sum_{i=1}^{N} \gamma_i  <\gamma$.
 \qedb
\end{pf}

%%%%
%\begin{proof}
%	Let $\gamma = \sum_{i=2}^{N}\gamma_i$ with $\gamma_i >0$ for all $i=2,\ldots,N$, then inequality (\ref{initial_condition_n_1}) is equivalent to 
%\begin{equation}\label{dis_cost}
%\bar{x}_{i0}^{\top} P_i \bar{x}_{i0}  <\gamma_i,\quad i= 2, 3, \ldots, N.
%\end{equation}
%
%
%(if)
%	For given feedback gain $\bar{K}$, if there exist symmetric solutions $P_i$ to inequalities (\ref{are_ineq}) and (\ref{dis_cost}) for all $i = 2,\ldots,N$, respectively, then, by Theorem 1, we have $J_i(\bar{K}) <\gamma_i$ for all $i = 2,\ldots,N$.
%Recall that $J(\bar{K}) = \sum_{i=2}^{N}{J}_i(\bar{K}) $ and $\gamma = \sum_{i=2}^{N}\gamma_i$, consequently, we have $J(\bar{K}) <\gamma$.
%
%(only if) 
%Recall that $J(\bar{K}) = \sum_{i=2}^{N}{J}_i(\bar{K}) $ and $\gamma = \sum_{i=2}^{N}\gamma_i$. If  $J(\bar{K}) <\gamma$, then $J_i(\bar{K}) <\gamma_i$ for all $i = 2,\ldots,N$, respectively. Again by applying Theorem 1, there exist symmetric solutions $P_i, i=2,3,\ldots,N$ to inequalities (\ref{are_ineq}) and (\ref{initial_condition_n_1}), respectively.
%\end{proof}
%\begin{rem}
%	 Lemma 2 only provides sufficient conditions for suboptimality, the reason is that we split the cost $J(K)$ equally to $N-1$ decoupled systems. Other choice could be, for example,  putting on weights for different agents based on their node degrees. 
%\end{rem}
Lemma \ref{lem_1} establishes a necessary and sufficient condition for a given feedback gain $K$ to stabilize all systems (\ref{closed_loop}) and to satisfy, for given initial states of these systems, $J(K)<\gamma$.
However, Lemma \ref{lem_1} does not yet provide a design method for computing such $K$.
Therefore, in the following we will provide a method to find such $K$.
%{\color{red}
%\subsection{Controllers that satisfy one global cost condition}
%}
%In Lemma 2, we need to solve $N-1$ Riccati inequalities to find $N-1$ solutions  $P_i, i = 2, 3, \ldots, N$, respectively. Next, we will simplify our design procedure and compute the control law.
%, we will establish conditions such that only one Riccati inequality is involved to find  one $P$ that is solution of all inequalities (\ref{are_ineq}). Moreover, we will also show how to obtain the suboptimal controller.
%{\color{blue}
%In the next theorem, we show that such a solution $P$ exists. Moreover, the corresponding feedback gain $K$ can be computed.
%}
%The following theorem yields a suboptimal control law for multi-agent system (\ref{net_sys}) that only one Riccati inequality is involved to compute the control gain. 

\begin{lem} \label{Main1}
	Consider the leader-follower system \eqref{sys_leader} and \eqref{sys_follower} with associated cost functional \eqref{cost}. Let $x_{r0}$ be the given initial state of the leader and $x_{i0}$, $i=1,2\ldots,N$ be the given initial states of the followers, respectively. Let $\gamma>0$ be a given tolerance.  
	Let $c$ be any real number such that  $0 < c < \frac{2}{\lambda_N}$. We distinguish two cases:
	\begin{enumerate}[(a)]
		\item \label{case1_lem}
		if 
		\begin{equation}\label{c1}
		\frac{2}{\lambda_1+\lambda_N} \leq c < \frac{2}{\lambda_N},
		\end{equation} 
		then there exists $P>0 $ satisfying 
		%\todoiny{?}{there exists a positive semi-definite matrix $P$ satisfying the Riccati inequality}
		\begin{equation} \label{one_are1}
		A^{\top}P + PA +(c^2\lambda_N^2-2c\lambda_N)PBR^{-1}B^{\top}P +\lambda_N Q < 0.
		\end{equation}
		\item \label{case2_lem}
		if 
		\begin{equation}\label{c2}
		0 < c < \frac{2}{\lambda_1+\lambda_N},
		\end{equation} 
		then there exists $P>0 $ satisfying 
		\begin{equation} \label{one_are2}
		A^{\top}P + PA +(c^2\lambda_1^2-2c\lambda_1)PBR^{-1}B^{\top}P +\lambda_N Q < 0.
		\end{equation}
	\end{enumerate}
	In both cases, if in addition $P$ satisfies
	\begin{equation}\label{p_N1}
	\sum_{i=1}^{N} (x_{i0} -x_{r0})^{\top} P (x_{i0} -x_{r0})< \gamma,
	\end{equation}
	then the distributed control law  \eqref{control} with $K := -cR^{-1}B^{\top}P$ achieves tracking consensus for the controlled leader-follower system  \eqref{sys_leader} and \eqref{sys_follower}, and with the initial states $x_{r0}$ and $x_{i0}$ we have  $J(K) <\gamma$.
\end{lem}

\begin{pf}
	We will only give the proof for case \ref{case1_lem} above. Using the upper and lower bounds on $c$ given by (\ref{c1}), it can be verified that $c^2\lambda_i^2-2c\lambda_i \leq c^2\lambda_N^2-2c\lambda_N <0$ for  $i=1,2,\ldots,N$. 
	It is then easily seen that \eqref{one_are1} has many positive definite solutions.
	Since also $\lambda_i \leq \lambda_N$, any such solution $P$ is a solution to the $N-1$ Riccati inequalities
	\begin{align}\label{n_are1}
	A^{\top}P + PA +(c^2\lambda_i^2-2c\lambda_i)PBR^{-1}B^{\top}P +\lambda_i Q &< 0, \nonumber \\ 
	i=1,2,\ldots,N.&
	\end{align}
	Equivalently, $P$ also satisfies
	%\todoiny{?}{$p$ also satisfies}
	the Lyapunov inequalities  
	\begin{equation}\label{lya_ineq}
	\begin{aligned}
	&(A-c\lambda_i BR^{-1}B^{\top}P)^{\top}P + P(A- c\lambda_iBR^{-1}B^{\top}P)  \\
	&\quad + \lambda_i Q + c^2\lambda_i^2 PBR^{-1}B^{\top}P < 0,\  i=1,2,\ldots,N.
	\end{aligned}
	\end{equation}
Next, by substituting $\bar{e} =(U^{\top}\otimes I_n)e$ into \eqref{initial_condition_n_1} we have 
$\sum_{i=1}^{N} e_{i0}^{\top} P e_{i0}  <\gamma$, which is equal to
 \eqref{p_N1}.
	
Next, taking $P_i = P$ for $i = 1,2,\ldots,N$ and $K := -cR^{-1}B^{\top}P$ in inequalities (\ref{are_ineq}) and (\ref{initial_condition_n_1}) immediately gives us inequalities (\ref{lya_ineq}) and (\ref{p_N1}).
Then it follows from Lemma \ref{lem_1} that all systems (\ref{closed_loop}) are stable and $J(K)<\gamma$. Subsequently,  the  controlled leader-follower system  reaches tracking consensus and $J(K)<\gamma$. \qedb
\end{pf}

We will now apply Lemma \ref{Main1} to establish a solution to Problem \ref{prob1}. The next  theorem provides a condition under which, for given radius $r$ and upper bound $\gamma$, suboptimal distributed diffusive control laws exist, and explains how these can be computed.
%{\color{red}
%	\begin{rem}
%		In Theorem \ref{Main1}, the controller design procedure is the following:
%		\begin{enumerate}
%			\item choose $c$ satisfying  \eqref{c1}, compute a solution $P$ of \eqref{one_are1} and obtain the corresponding control gain $K$;
%			\item for given initial state $x_0$ and upper bound $\gamma$ for the cost, check whether \eqref{p_N1} is satisfied;
%			\item if  \eqref{p_N1} is satisfied, then the control law $L \otimes  K$ is suboptimal. 
%			%
%			Otherwise, one needs to make the differences of the initial conditions between different agents relatively ``smaller'' in the sense as explained in Remark \ref{rem_main1}, or one needs to increse the desired tolerance $\gamma$.
%		\end{enumerate}
%	\end{rem}
%}
%
\begin{thm} \label{Realmain1}
	Consider the leader-follower system  \eqref{sys_leader} and \eqref{sys_follower}  with associated cost functional \eqref{cost}. Let $r >0$ be a given radius and let $\gamma >0$ be an a priori given upper bound for the cost.
	Let $c$ be any real number such that  $0 < c < \frac{2}{\lambda_N}$. We distinguish two cases:
	\begin{enumerate}[(a)]
		\item \label{case1_thm}
		if 
		\begin{equation}\label{c1new}
		\frac{2}{\lambda_1+\lambda_N} \leq c < \frac{2}{\lambda_N},
		\end{equation} 
		then there exists $P>0$ satisfying
		%\todoiny{?}{there exists a positive semi-definite matrix $P$ satisfying the Riccati inequality}
		\begin{equation} \label{one_are1new}
		A^{\top}P + PA +(c^2\lambda_N^2-2c\lambda_N)PBR^{-1}B^{\top}P +\lambda_N Q < 0.
		\end{equation}
		\item \label{case2_thm}
		if 
		\begin{equation}\label{c2new}
		0 < c < \frac{2}{\lambda_1+\lambda_N},
		\end{equation} 
		then there exists $P>0 $ satisfying 
		\begin{equation} \label{one_are2new}
		A^{\top}P + PA +(c^2\lambda_1^2-2c\lambda_2)PBR^{-1}B^{\top}P +\lambda_N Q < 0.
		\end{equation}
	\end{enumerate}
	In both cases, if in addition $P$ satisfies
	\begin{equation} \label{ineq}
	P < \frac{\gamma}{r^2} I,
	\end{equation}
	then the distributed control law \eqref{control} with $K := -cR^{-1}B^{\top}P$  achieves tracking consensus for the controlled leader-follower system  \eqref{sys_leader} and \eqref{sys_follower}  and $J(K) <\gamma$ for all initial states  $x_{r0}$ and $x_0$  satisfying 
	\begin{equation}\label{initial}
			x_{0} - \mathbf{1}_N \otimes   x_{r0} \in B(r).
	\end{equation}
\end{thm}
\begin{pf}
	Again, we only give proof for case \ref{case1_thm}. Let $P>0$ satisfy \eqref{one_are1new} and \eqref{ineq}. Next, we will show that if the initial states  $x_{r0}$ and $x_0$  satisfy  $x_{0} - \mathbf{1}_N \otimes   x_{r0} \in B(r)$, then \eqref{p_N1} holds.
	Indeed, if $\|x_{0} - \mathbf{1}_N \otimes   x_{r0} \| \leq r$, then
	\begin{equation*}
	\begin{aligned}
		&\sum_{i=1}^{N} (x_{i0} -x_{r0})^{\top} P (x_{i0} -x_{r0}) \\
	=	&  (x_{0} - \mathbf{1}_N\otimes x_{r0})^\top (I \otimes P) (x_{0} - \mathbf{1}_N \otimes x_{r0}) \\
	<	& \frac{\gamma}{r^2} \|x_{0} - \mathbf{1}_N \otimes   x_{r0} \|^2 \leq \gamma.
	\end{aligned}
	\end{equation*}
	It then follows from Lemma \ref{Main1} that the controlled leader-follower system   \eqref{sys_leader} and \eqref{sys_follower}  reaches tracking consensus with the given $K$  and $J(K) <\gamma$ for all initial states  $x_{r0}$ and $x_0$  satisfying \eqref{initial}.
\qedb
\end{pf}

\begin{rem} \label{rem_main1}
	Theorem \ref{Realmain1} states that after choosing $c$ satisfying the inequality (\ref{c1new}) for case \ref{case1_thm} and finding  $P>0$ satisfying (\ref{one_are1new}) and \eqref{ineq}, the distributed control law with local gain $K = -c R^{-1}B^{\top}P$ is $\gamma$-suboptimal for all initial states of the leader-follower system satisfying the condition \eqref{initial}. 
	According to \eqref{ineq}, the smaller the solution $P$ of \eqref{one_are1new}, the smaller the quotient $\frac{\gamma}{r^2}$ is allowed to be, leading to a smaller upper bound and a larger radius.
	The question then arises: how should we choose the parameter $c$ in \eqref{c1new} so that the Riccati inequality \eqref{one_are1new} allows a positive definite solution that is as small as possible.
	In fact, one can find a positive definite solution $P(c,\epsilon)$ to (\ref{one_are1new}) by solving the Riccati equation
	\begin{equation}\label{are_rem}
	A^{\top}P + PA -PB\bar{R}(c)^{-1}B^{\top}P +\bar{Q}(\epsilon)= 0
	\end{equation} 
	with  $\bar{R}(c) = \frac{1}{-c^2\lambda_N^2+2c\lambda_N}R$ and $\bar{Q}(\epsilon) =\lambda_N Q +\epsilon I_n $	where  $c$ is chosen as in (\ref{c1new}) and  $\epsilon >0$.
	If $c_1$ and $c_2$ as in (\ref{c1new}) satisfy $c_1 \leq c_2$, then we have
	$\bar{R}(c_1) \leq \bar{R}(c_2)$,
	so, clearly, $P(c_1,\epsilon) \leq P(c_2,\epsilon)$.
	Similarly, if $0 <  \epsilon_1 \leq \epsilon_2$,  we immediately have  $\bar{Q}(\epsilon_1) \leq \bar{Q}(\epsilon_2)$. Again, it follows that $P(c,\epsilon_1) \leq P(c,\epsilon_2)$.
	Therefore, if we choose $\epsilon>0$ very close to $0$ and $c = \frac{2}{\lambda_1+\lambda_N}$,  we find the `best' solution to the Riccati inequality (\ref{one_are1new}) in the sense explained above.

	Likewise, if $c$ satisfies \eqref{c2new} corresponding to case \ref{case2_thm}, it can be shown that if we choose $\epsilon>0$ very close to $0$ and $c>0$ very close to $\frac{2}{\lambda_1+\lambda_N}$,  we find the `best' solution to the Riccati inequality (\ref{one_are2new}) in the sense explained above.
\end{rem}

%\begin{rem}
%Give a remark here to say something on why would we consider the initial states in a ball.
%\end{rem}
%\begin{rem}
%	Is there result to estimate the eigenvalues of $\Gamma$?
%\end{rem}
%\begin{rem}
%May be also result on weighted graph?
%\end{rem}

\section{Simulation Example}\label{sec_simulation}
%{\color{red}
In this section, we will use a numerical example borrowed from \cite{Nguyen2015} to illustrate the design method for the suboptimal distributed control laws given in Theorem \ref{Realmain1}.

\begin{figure}[t]
	\centering
	\begin{tikzpicture}[scale=0.7]
	\tikzset{VertexStyle1/.style = {shape = circle,
			color=black,
			fill=white!93!black,
			minimum size=0.5cm,
			text = black,
			inner sep = 2pt,
			outer sep = 1pt,
			minimum size = 0.55cm},
		VertexStyle2/.style = {shape = circle,
			color=black,
			fill=black!53!white,
			minimum size=0.5cm,
			text = white,
			inner sep = 2pt,
			outer sep = 1pt,
			minimum size = 0.55cm}
	}
	\node[VertexStyle2,draw](0) at (-4,0) {$\bf r$};
	\node[VertexStyle1,draw](2) at (-1,0) {$\bf 2$};
	\node[VertexStyle1,draw](3) at (1.2,2) {$\bf 3$};
	\node[VertexStyle1,draw](4) at (4,1.3) {$\bf 4$};
	\node[VertexStyle1,draw](5) at (4,-1.3) {$\bf 5$};
	\node[VertexStyle1,draw](1) at (1.2,-2) {$\bf 1$};
	\Edge[ style = {->,> = latex',pos = 0.2},color=black, labelstyle={inner sep=0pt}](0)(2);
	\Edge[ style = {-,> = latex',pos = 0.2},color=black, labelstyle={inner sep=0pt}](1)(2);
	\Edge[ style = {-,> = latex',pos = 0.2},color=black, labelstyle={inner sep=0pt}](2)(3);
	\Edge[ style = {-,> = latex',pos = 0.2},color=black, labelstyle={inner sep=0pt}](3)(4);
	\Edge[ style = {-,> = latex',pos = 0.2},color=black, labelstyle={inner sep=0pt}](4)(5);
	\Edge[ style = {-,> = latex',pos = 0.2},color=black, labelstyle={inner sep=0pt}](5)(1);
	\end{tikzpicture}
	\caption{The underlying graph of the communication between the leader and the followers.}
	\label{graph}
\end{figure}
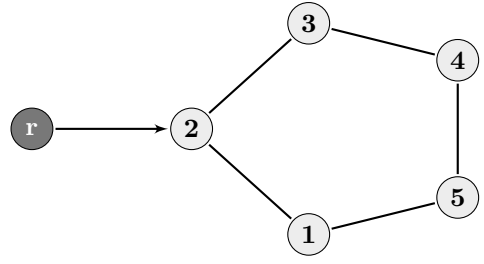

Consider a leader-follower multi-agent system, consisting of one  leader and five followers.
The dynamics of the leader is given by
\begin{equation*}
\dot{x}_r(t) = Ax_r(t),\quad x_r(0) = x_{r0},
\end{equation*}
and the dynamics of the followers are identical and  represented by
\begin{equation*}
\dot{x}_i(t) = Ax_i(t) + Bu_i(t),\quad x_i(0) = x_{i0}, \quad i=1, 2,\ldots,5
\end{equation*}
where 
\begin{equation*}
A = 
\begin{pmatrix}
0 & 1 \\
-1 & 0
\end{pmatrix},\quad 
B = 
\begin{pmatrix}
0 \\ 1
\end{pmatrix}
.
\end{equation*}
The pair $(A,B)$ is stabilizable.
Assume the underlying graph representing the communication between the leader and the followers is given as in Figure \ref{graph}.
The  graph representing the communication between the followers is then the undirected cycle graph with the Laplacian matrix 
\begin{equation*}
L = 
\begin{pmatrix}
2 & -1 & 0 & 0 & -1 \\
-1 & 2& -1 & 0 & 0 \\
0& -1 & 2 & -1 & 0 \\
0 & 0 & -1 & 2 & -1 \\
-1 & 0 & 0 & -1 & 2
\end{pmatrix}.
\end{equation*}
%
%\noindent

Since the leader shares its state information only with follower $2$,  it follows from Lemma \ref{Gamma} that the associated diagonal matrix $G =\textnormal{diag}(g_1,g_2,\ldots,g_5) = \textnormal{diag}(0,1,0,0,0).$
Furthermore, we consider the cost functional
\begin{equation*} 
\begin{aligned}
J(u) &= \int_{0}^{\infty}\frac{1}{2}\sum_{i=1}^{5}\sum_{j=1}^5 a_{ij}(x_i-x_j)^{\top} Q (x_i-x_j) \\
& \quad + \sum_{i=1}^{5}  g_{i}(x_i-x_r)^{\top} Q (x_i-x_r) + \sum_{i=1}^{N}u_i^{\top}R u_i \ dt
\end{aligned}
\end{equation*}
with 
\begin{equation*}
Q = 
\begin{pmatrix}
	2 & 0 \\
	0 & 1
\end{pmatrix},\quad
R = 1.
\end{equation*}
Let the desired tolerance for the cost functional be $\gamma = 20$.
Our aim is then to design a control law of the form 
\begin{equation}\label{control_simulation} 
u_i (t)=  K\sum_{j=1}^{5}a_{ij} (x_i(t) - x_j (t)) +K g_i (x_i(t) - x_r(t))
\end{equation} 
such that the controlled leader-follower system reaches tracking consensus and the associated cost satisfies $J(K) < 20$ for all initial states $x_0$ and $x_{r0}$ satisfying the condition $\|x_{0} - \mathbf{1}_5\otimes   x_{r0}\| \leq r$ with radius $r$ to be specified later.

In this simulation example, we will use the design method of case \ref{case1_thm} in Theorem \ref{Realmain1}.
For $\Gamma = L+G$ the smallest and largest eigenvalues are $\lambda_1 =0.1392$ and $\lambda_5 =4.1149$, respectively.
We first compute a solution $P>0$ to  \eqref{one_are1new} by solving
\begin{equation}\label{are}
	A^{\top}P + PA +(c^2\lambda_5^2-2c\lambda_5)PBR^{-1}B^{\top}P +\lambda_5 Q +\epsilon I_2 = 0
\end{equation}
with $\epsilon$ sufficiently small as mentioned in Remark \ref{rem_main1}. Here we choose $\epsilon =0.01$.
In addition, we choose $c = \frac{2}{\lambda_1+\lambda_5} =0.4701$, which is the `best' choice as mentioned in Remark \ref{rem_main1}.
Then by solving \eqref{are} using Matlab, we compute 
\begin{equation*}
P = 
\begin{pmatrix}
   13.2553   & 3.3886\\
3.3886  &  9.2760
\end{pmatrix}.
\end{equation*}
Correspondingly, the control gain is equal to
$K = \begin{pmatrix}
    1.5931   & 4.3610
\end{pmatrix}$.
We now compute the radius $r$ of a ball $B(r)$ of initial states for which the distributed control law \eqref{control_simulation} is suboptimal, i.e. $J(K) < 20$. 
We compute that the largest eigenvalue of $P$ is equal to $15.1952$. Hence for every radius $r$ such that $\frac{20}{r^2} > 15.1952$ the inequality \eqref{ineq} holds.  Thus, the distributed controller with local  gain $K$  is suboptimal for all $x_{r0}$ and $x_0$ satisfying $\|x_{0} - \mathbf{1}_5\otimes   x_{r0}\| \leq r$ with $r < 1.1473
$.

\begin{figure}[t]
%	\centering
	\includegraphics[width=\columnwidth]{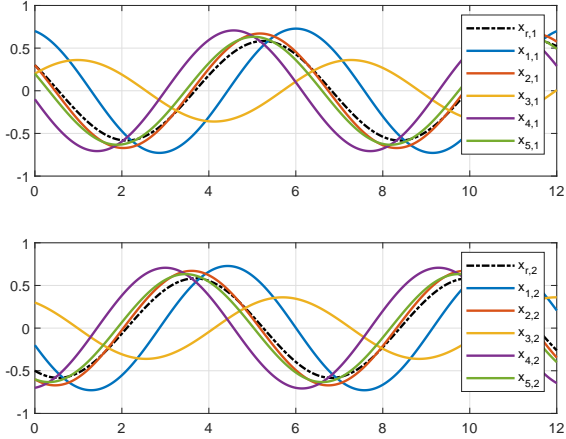}
	\caption{Plots of the states $x_{r1}$ and $x^1 = (x_{1,1},\ldots, x_{5,1})$ (upper plot) and  $x_{r2}$ and $x^2 = (x_{1,2},\ldots, x_{5,2})$ (lower plot) of the six decoupled  local agents without control} \label{decoupled}
\end{figure}

\begin{figure}[t]
	\centering
	\includegraphics[width=\columnwidth]{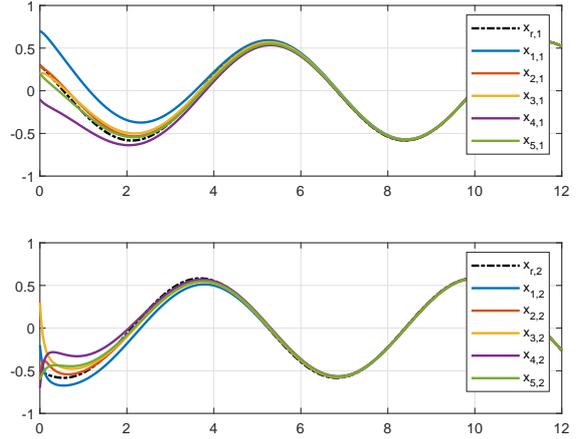}
	\caption{Plots of the states $x_{r1}$ and $x^1 = (x_{1,1},\ldots, x_{5,1})$ (upper plot) and $x_{r2}$ and $x^2 = (x_{1,2},\ldots, x_{5,2})$ (lower plot) of the controlled leader-follower system} \label{consensus}
\end{figure}

As an example, the following initial states of the agents satisfy this norm bound: 
$x_{r0}^{\top} = 
\begin{pmatrix} 
0.3 & -0.5
\end{pmatrix}$,
$x_{10}^{\top} = 
\begin{pmatrix} 
 0.7 & -0.2
\end{pmatrix}$,
$x_{20}^{\top} = 
\begin{pmatrix} 
0.3 & -0.6
\end{pmatrix}$,
$x_{30}^{\top} = 
\begin{pmatrix} 
0.2 & 0.3
\end{pmatrix}$,
$x_{40}^{\top} = 
\begin{pmatrix} 
-0.1  & -0.7
\end{pmatrix}$,
$x_{50}^{\top} = 
\begin{pmatrix} 
0.2 & -0.6
\end{pmatrix}$.
The plots of the state of the six local agents without control are shown in Figure \ref{decoupled}.
Figure \ref{consensus} shows that the controlled leader-follower system reaches tracking consensus.

\section{Conclusion}\label{sec_conclusion}
In this paper, we have studied the distributed linear quadratic tracking control problem for leader-follower multi-agent systems. We have considered leader-follower systems consisting of one autonomous leader and $N$ followers, together with an associated  global cost functional. 
We assume that the leader shares its state information with at least one of the followers and the underlying graph connecting the followers is a connected simple undirected graph.
For this type of leader-follower systems, we have provided a design method to compute distributed suboptimal control laws such that the controlled network reaches tracking consensus and the associated cost is smaller than a given tolerance for all initial states bounded in norm by a given radius.
The computation of the local gain involves the solution of a single Riccati inequality, whose dimension is equal to the dimension of the agent dynamics, and also involves the largest and smallest eigenvalue of a positive definite matrix capturing the underlying graph structure.

%A possible direction for future work is to consider the cases that leader is not autonomous but with bounded input signal (\cite{bibid}), another interesting direction is to consider the case with multiple leaders, which is called containment control (\cite{bibid}).
%Another possible direction is to consider heterogeneous agent dynamics instead of homogeneous agent dynamics. 
\balance

%\newpage
%\bibliography{ifacconf}             % bib file to produce the bibliography
%\bibliographystyle{IEEEtran}
\bibliography{suboptimal}                         % with bibtex (preferred)
                                                   
%\begin{thebibliography}{xx}  % you can also add the bibliography by hand

%\bibitem[Able(1956)]{Abl:56}
%B.C. Able.
%\newblock Nucleic acid content of microscope.
%\newblock \emph{Nature}, 135:\penalty0 7--9, 1956.

%\bibitem[Able et~al.(1954)Able, Tagg, and Rush]{AbTaRu:54}
%B.C. Able, R.A. Tagg, and M.~Rush.
%\newblock Enzyme-catalyzed cellular transanimations.
%\newblock In A.F. Round, editor, \emph{Advances in Enzymology}, volume~2, pages
%  125--247. Academic Press, New York, 3rd edition, 1954.

%\bibitem[Keohane(1958)]{Keo:58}
%\newblock \emph{Power and Interdependence: World Politics in Transitions}.
%\newblock Little, Brown \& Co., Boston, 1958.

%\bibitem[Powers(1985)]{Pow:85}
%\newblock Is there a way out?
%\newblock \emph{Harpers}, pages 35--47, June 1985.

%\bibitem[Soukhanov(1992)]{Heritage:92}
%\newblock \emph{{The American Heritage. Dictionary of the American Language}}.
%\newblock Houghton Mifflin Company, 1992.

%\end{thebibliography}
%
%\appendix
%\section{A summary of Latin grammar}    % Each appendix must have a short title.
%\section{Some Latin vocabulary}              % Sections and subsections are supported  
                                                                         % in the appendices.
\end{document}